\def \Reg {{\rm Reg}}
\def \Qed {\rule{2.5mm}{3mm}}
\def\ZZ{{\hbox{\sf Z\kern-.43emZ}}}
\def \H {{\rm H}}
\def \ot {\otimes}
\def \Z {{\rm Z}}
\def \B {{\rm B}}
\def \kb {k^{\bullet}}

\def \De {\Delta}
\def \ep {\varepsilon}
\def \id {{\rm id}}

\def \Pp {{\rm P}}
\def \Hom {{\rm Hom}}
\def \res {{\rm res}}
\def \Hm {{\rm H_m^2}}
\def \Hc {{\rm H_c^2}}
\def \io {\iota}
\def \m {{\rm m}}
\def \Map {{\rm Map}}
\def \om {\omega}

\def \al {\alpha}

\def \Alg {{\rm Alg}}
\def \g   {{\bf g}}
\def \h  {{\bf h}}
\def \Regm {{\rm Reg}_{\rm meas}}
\def \Alg {{\rm Alg}}
\def \Lie {{\rm Lie}}
\def \Vec {{\rm Vect}}

\def \U {{\cal U}}
\def \M {{\rm M}}
\def \th {\theta}
\def \ol {\overline}
\def \de {\delta}
\def \Coalg {{\rm Coalg}}
\def \rt {\rtimes}

\documentclass[10 pt]{article}
\title{On the Cohomology of a Smash Product of Hopf Algebras}
\author{Mitja Mastnak \\ Department of Mathematics and
Statistics\\Dalhousie University\\ Halifax, N.S., Canada, B3H
3J5\\E-mail: mastnak@mscs.dal.ca }

\usepackage{amssymb}

\newtheorem{Theorem}{Theorem}[section]
\newtheorem{Lemma}[Theorem]{Lemma}
\newtheorem{Proposition}[Theorem]{Proposition}
\newtheorem{Definition}[Theorem]{Definition}
\newtheorem{Corollary}[Theorem]{Corollary}

\newtheorem{Example}[Theorem]{Example}

\begin{document}

\tolerance=1000

\maketitle

\begin{abstract}
A five term sequence for the low degree cohomology of a smash product of
(cocommutative) Hopf algebras is obtained, generalizing that of Tahara for a
semi-direct product of groups
\end{abstract}

\setcounter{section}{-1}

\section{Introduction}

Our aim is to obtain a five term sequence for the cohomology of a
smash product of (cocommutative) Hopf algebras, which generalizes
the Tahara sequence for a semi-direct product of groups [Ta].

If the finite group $G$ is a semi direct product of subgroups $N$
and $T$ and $A$ is an abelian group, which is also a trivial $\ZZ
G$ module then Tahara's exact sequence is
\begin{eqnarray*}
1 &\to& \H^1(T,\H^1(N,A))\to \tilde{\rm
H}^2(G,A)\stackrel{\res}{\to}\H^2(N,A)^T\stackrel{d_2}{\to} \\
&{\to}&\H^2(T,\H^1(N,A))\to \tilde{\rm H}^3(G,A).
\end{eqnarray*}

In our setting $N$ and $T$ are cocommutative Hopf algebras (hence
$G$ a smash product of $N$ and $T$) and $A$ is replaced by a
commutative algebra which is also a trivial $G$-module. The
cohomology in question is the Sweedler cohomology of an algebra
over a Hopf algebra. The Tahara sequence then generalizes to
\begin{eqnarray*}
0&\to& \H^1_{\rm meas}(T,\Hom(N,A))\stackrel{\io}{\to} \tilde{\rm
H}^2(G,A)\stackrel{\res}{\to}
\H^2(N,A)^T\\
&\stackrel{d}{\to}& \H^2_{\rm meas}(T,\Hom(N,A))\stackrel{j}{\to}
\tilde{\rm H}^3(G,A),
\end{eqnarray*}
where $\H^2(N,A)^T$ is the Hopf algebra analogue for the $T$-stable part
of cohomology and $\H^i_{\rm meas}$ denotes the measuring cohomology
(generated by the regular maps that make $T$ measure $N$ to $A$).
We can identify $\H^2_{\rm
meas}(T,\Hom(N,A))$ with $\Hm(T,M(N,A))$, the
multiplication part of cohomology of a Singer pair [M], and $\H^1_{\rm
meas}(T,\Hom(N,A))$ with a certain quotient of $\H^1(T,M(N,A))$ ($M(N,A)$ is
the universal measuring Hopf algebra).
This sequence is a generalization of the
Tahara sequence in the sense that $\H^i(kH,A)\simeq \H^i(H,{\cal
U}(A))$, where $kH$ is a group algebra and ${\cal U}(A)$ denotes
the abelian group of units of $A$.

\section{Preliminaries}

All vector spaces are over a fixed ground field $k$.

\subsection{Smash Product of Hopf Algebras}\label{s1}

Let $T$ and $N$ be cocommutative Hopf algebras and let $T$ act on
$N$ via $@\colon T\ot N\to N$, we shall abbreviate $@(t\ot
n)=t(n)$, so that $N$  becomes a $T$-module bialgebra. The smash
product $H=N\rtimes T$ (also called crossed product), a special case
of a bicrossed product (or bismash product) arising from a
matched pair (see [Kas, IX.2]), is the tensor product coalgebra
$N\ot T$ (to avoid confusion we shall write $n\# t$
for $n\ot t$) with multiplication given by:
\begin{eqnarray*}
(n\# t)(n'\# t')&=&nt(n')\# tt'.
\end{eqnarray*}
We identify $N$ and $T$ with $N\# 1$ and $1\# T$ respectively  as
Hopf subalgebras of $H$. In this manner we have the equality $n\#
t=nt$ for $n\in N$ and $t\in T$.

\subsection{The Sweedler cohomology}\label{s2}

The Sweedler cohomology is a cohomology theory for algebras that
are modules over a given Hopf algebra [Sw1]. Let $H$ be a
cocommutative Hopf algebra, and let $A$ be a commutative algebra,
that is also a (right) $H$-module. Although it is more usual to consider a left action, we
work with a right one in this paper, for the sake of convenience (one can make a left action out
of a right one by using the antipode).
If the $H$-module structure is given
by a right action $\Psi\colon A\ot H\to A$, $a\ot h\mapsto a^h$, then
the standard (normalized) complex for computing the Sweedler cohomology is
as follows:
$$ \ldots \to \Reg(H^{\ot q-1},A)\stackrel{\delta^{q-1}}{\to}
\Reg(H^{\ot q},A)\to \ldots ,$$ where
the differential is
given by $\delta^{q-1}(f)=(\ep\ot f)*f^{-1}(m\ot\id\ot\ldots \ot\id)*
\ldots *f^{\pm}(\id\ot\ldots\ot\id\ot\m)*\Psi(f^{\mp}\ot\id)$.
Here $\Reg(H^{\ot q},
A)$ denotes the abelian group of convolution invertible
linear maps $f\colon H^{\ot q}\to A$ with the property that
$f(h_1\ot \ldots h_q)=\ep(h_1)\ldots \ep(h_q)1_A$ whenever some
$h_i\in k$. The cocycles and coboundaries for the degree $1$ and
degree $2$ cohomology groups are described as follows:
\begin{eqnarray*}
\Z^1(H,A)&=&\{f\in\Reg(H,A)|f(x y) = f(y_1)f(x)^{y_2}\}, \\
\B^1(H,A)&=&\{f\in\Reg(H,A)|\exists a\in {\cal U}(A), \;{\rm
s.t.}\;
f(x)=a(a^{-1})^x\}, \\
\Z^2(H,A)&=&\{f\in\Reg(H\ot H,A)|f(x_1\ot y_1)^{z_1}
f(x_2y_2\ot z_2)=\\
&&=f(y_1\ot z_1)f(x \ot y_2 z_2) \}, \\
\B^2(H,A)&=&\{f\in\Reg(H\ot H,A)|\exists t\in \Reg(H,A),\;{\rm
s.t.}\\
&& f(x\ot y)=t(y_1)t^{-1}(x_1y_2)t(x_2)^{y_3}\}.
\end{eqnarray*}

The second Sweedler cohomology group
$\H^2(H,A)=\Z^2(H,A)/\B^2(H,A)$ classifies equivalence classes of
cleft comodule algebra extensions $A\stackrel{\io}{\to}
C\stackrel{\pi}{\to} H$ [Sw1, Mo]. These are sequences of
$H$-comodule algebra maps with the property that $A=C^{co H}=
\{c\in C| \rho(c)=c\ot 1\}$(here $\rho\colon C\to C\ot H$ denotes
the $H$-comodule structure on $C$) and that there exists a
convolution invertible $H$-comodule map $\chi\colon H\to C$ (which
is automatically a splitting of $\pi$). The isomorphism between
extensions and cohomology is induced by $\chi\mapsto f_\chi$,
where $f_\chi\in\Z^2(H,A)$ is the Sweedler cocycle given by
$f_\chi(h\ot h')=\chi(h_1)\chi(h'_1)\chi^{-1}(h_2h'_2).$

\subsection{Cohomology of a Singer pair}

A Singer pair (in literature also the term matched pair is used [Ho, Si], the term Singer pair was
introduced by A. Masuoka [Ma]) consists of a pair of Hopf algebras $(B,A)$ together with an
action and a coaction that are compatible in some sense. For the sake of convenience we shall
consider right Singer pairs, i.e.
the ones given by a right action $\mu\colon A\ot B\to A$ and a right coaction
$\rho\colon B\to A\ot B$.

In the special case when the coaction is trivial,
a (right) Singer pair is given by a right action $\mu\colon A\ot B\to A$ which makes
$A$ into a $B$-module bialgebra. If the Hopf algebra $A$ is commutative and the
Hopf algebra $B$ is cocommutative then the Singer pair is called abelian.
In this paper we consider (right) abelian Singer pairs
with the trivial coaction exclusively.

The cohomology of an abelian Singer pair $(B,A)$ is computed by the total complex of a
certain double complex [Ho]. Cocycles $\Z^n(B,A)$ (coboundaries $\B^n(B,A)$) are
$n$-tuples of linear maps $(f_i)_{i=1}^n$,
$f_i\colon  B^{\ot i}\to A^{\ot n+1-i}$ that satisfy cocycle (coboundary) conditions.

In order to compare measuring cohomology (to be defined in Section \ref{s21})
to cohomology of Singer pairs, we need to
introduce subgroups $Z^n_i(B,A)$ of $\Z^n(B,A)$, that are spanned
by $n$-tuples in which the $f_j$'s are trivial for $j\not= i$ and
subgroups $\B^n_i=\Z^n_i\cap \B^n\le \B^n$. These give rise to
subgroups of cohomology groups $\H^n_i=\Z^n_i/\B^n_i\simeq
(\Z^n_i+\B^n)/\B^n\le \Z^n/\B^n=\H^n$. If $n=2$ then
$\H^2_1(B,A)=\Hc(B,A)$ is the comultiplication
part and $\H^2_2(B,A)=\Hm(B,A)$ is the multiplication part of
cohomology [M].

We are mainly interested in subgroups $\H_n^n(B,A)$ of
$\H^n(B,A)$. If $n\ge 2$ they can be computed by a subcomplex of
the complex in Section \ref{s2}. This is easily seen by comparing
that complex to the complex for computing cohomology of a Singer
pair [Ho, Ma] (recall that the coaction $\rho\colon B\ot A\to B$
is trivial).
\begin{Proposition}\label{desc}
If $n\ge 2$ then we have canonical isomorphisms
$$\Z^n_n(B,A)\simeq \{f\in \Coalg(B,A)|\de^n f=\ep\},$$
where $\de^n f\colon \Reg(B^{\ot n},A)
\to \Reg(B^{\ot n+1},A)$ is the Sweedler differential (see Section \ref{s2}), and
$$\B^n_n(B,A)\simeq\{f\colon B^{\ot n}\to A|\exists g\in \Coalg(B^{\ot n-1},A),\; s.t.\;
f=\de^{n-1} g\}.$$
\end{Proposition}

The second cohomology group of a Singer pair $(B,A,\mu,\rho)$ classifies equivalence classes of
those Hopf algebra extensions that induce the same Singer pair [Ho]. Hopf algebra extensions are
sequences of Hopf algebra maps $A\stackrel{\io}{\to} C\stackrel{\pi}{\to} B$, such that
$A=C^{co B}=\{c\in C|(\id\ot \pi)\De c=c\ot 1\}$ and there exists a $B$-comodule map $\chi\colon
B\to A$ or equivalently $B\simeq C_A=C/\ep_A^{-1}(0)C$ and there exists an $A$-module map
$\xi\colon C\to A$. The group $\H^2_2(B,A)=\Hm(B,A)$ classifies cocentral Hopf algebra extensions
$A\to C\to B$ for which there exists an $A$-module coalgebra map $\xi\colon C\to A$.

{\bf\noindent Remark. }{\sl We can make a left Singer pair out of the right one with the
use of antipodes. The cohomology groups for such pairs would remain the same (isomorphic).}

\subsection{Universal measuring coalgebra}
Here we introduce the (contravariant) universal coalgebra functor
$\M(\_,A)\colon  {\bf Alg}^{op}\linebreak \to {\bf Coalg}$ from the category of
algebras to the category of coalgebras,
that is adjoint to $\Hom(\_,A)\colon {\bf Coalg}\to {\bf Alg}^{op}.$ For
more details, we refer to [Sw2, Chapter VII] and
[GP].

Let $A$, $B$, be algebras, $C$ a coalgebra.
\begin{Proposition}[Sw2, 7.0.1]
A map $\psi\colon C\ot B\to A$ corresponds to an algebra map $\rho\colon
B\to\Hom(C,A)$, $\rho(b)(c)=\psi(c\ot b)$ if and only if
\begin{enumerate}
\item
$\psi(c\ot bb')=\psi(c_1\ot b)\psi(c_2\ot b'),$
\item
$\psi(c\ot 1_B)=\ep(c)1_A$
\end{enumerate}
\end{Proposition}
If the equivalent conditions from the proposition above are
satisfied, we say that $(\psi,C)$ {\bf measures} $B$ to $A$.

Given algebras $B$ and $A$ there is a measuring $(\theta,\M(B,A))$
with the following universal property.
\begin{Theorem}[Sw2, 7.0.4]
The universal measuring $\theta\colon \M(B,A)\ot B\to A$ has the
following universal property:

for any measuring $f\colon C\ot B\to A$ there exists a unique
coalgebra map $\overline{f}\colon C\to \M(B,A)$, s.t.
$f=\theta(\overline{f}\ot\id)$.
\end{Theorem}

\section{The measuring cohomology}\label{s21}

Let $H=N\rtimes T$ (as in Section \ref{s1}) and let $T$ act on
$\Hom(N,A)$ via pre-composition: $\Psi\colon \Hom(N,A)\ot T\to
\Hom(N,A)$, $f\ot t\to f^t$, where $f^t(n)=f(t(n))$. Let
$\Reg_{\rm meas}(T^{\ot q},\Hom(N,A))$ denote the subgroup of
$\Reg(T^{\ot q},\Hom(N,A))$ consisting of the maps $f\colon T^{\ot
q}\to \Hom(N,A)$ that make $T^{\ot q}$ measure $N$ to $A$,
i.e. $f({\bf t})(nn')=\sum f({\bf
t}_1)(n)f({\bf t}_2)(n')$ and $f({\bf t})(1_N)=\ep({\bf t})1_A$ for
${\bf t}\in T^{\ot q}$ and $n,n'\in
N$. The differential $$\Reg(T^{\ot q-1},
\Hom(N,A))\stackrel{\delta^{q-1}}{\to} \Reg(T^{\ot
q},\Hom(N,A)),$$ described in Section \ref{s2} restricts to
$$\Reg_{\rm meas}(T^{\ot q-1},\Hom(N,A))\stackrel{\delta^{q-1}}{\to}
\Reg_{\rm meas}(T^{\ot q},\Hom(N,A)),$$ thus giving rise to a sub
complex of the complex given in Section \ref{s2}. We name the
cohomology it produces the "measuring cohomology" and denote it by
$\H_{\rm meas}^q(T,\Hom(N,A))$. We denote the measuring cocycles
and coboundaries by $\Z_{\rm meas}^q(T,\Hom(N,A))$ and $\B_{\rm
meas}^q(T,\Hom(N,A))$, respectively. If $q=1$ or $q=2$ these are as follows:
\begin{eqnarray*}
\Z^1_{\rm meas}&=&\{f\in \Reg_{\rm meas}(T,\Hom(N,A))|
f(ts)(n)=f(s_1)(n_1)f(t)(s_2(n_2))\}\\
\B^1_{\rm meas}&=&\{f\in \Reg_{\rm meas}(T,\Hom(N,A))|
\exists g\in \Hom(N,A)\;{\rm s.t.}\\
&&f(t)(n)=g(n_1)g^{-1}(t(n_2))\}\\
\Z^2_{\rm meas}&=&\{f\in \Reg_{\rm meas}(T\ot T,\Hom(N,A))|
f(t_1\ot s_1)(r_1(n_1))\\
&&\cdot f(t_2s_2\ot r_2)(n_2)=f(s_1\ot r_1)(n_1)f(t\ot s_2r_2)(n_2)\}\\
\B^2_{\rm meas}&=&\{f\in \Reg_{\rm meas}(T\ot T,\Hom(N,A))|
\exists g\in \Reg_{\rm meas}(T,\Hom(N,A))\; \\
&&{\rm s.t.}\;f(ts)(n)=f(s_1)(n_1)f(t)(s_2(n_2))\}
\end{eqnarray*}
We will often identify $f({\bf t})(n)=f({\bf t}\ot n)$ for
$f\in \Hom(T^{\ot q},\Hom(N,A))\simeq\Hom(T^{\ot q}\ot N,A).$

{\bf\noindent Remark. }{\sl Even though the notation $\H^i_{\rm
meas}(T,N,A)$ would probably be more precise (since $\H^i_{\rm
meas}(T,B)$ only makes sense if $B=\Hom(N,A)$ and it is also possible that
$\Hom(N,A)=\Hom(N',A')$ for different $N, N', A, A'$), we write
$\H^i_{\rm meas}(T,\Hom(N,A))$ to emphasize the fact that the
measuring cohomology is computed by a subcomplex of the standard
complex for computing Sweedler cohomology. }

\subsection{$T=kG$}
In the case $T=kG$ is a group algebra we can
identify
$$\Regm((kG)^{\ot i},\Hom(N,A))\simeq \Map(G^{\times i},\Alg(N,A)),$$
where $\Alg(N,A)$ denotes the group of algebra maps from $N$ to
$A$ (with convolution product) and $\Map(H,B)$ denotes the abelian
group of unital maps from $H$ to $B$ (with pointwise
multiplication). Note that this isomorphism induces an isomorphism
of complexes, i.e. preserves the differentials; hence we have:
\begin{Theorem}
There is an isomorphism
$$\H^i_{\rm meas}(kG,\Hom(N,A))\simeq
\H^i(G,\Alg(N,A)).$$
\end{Theorem}

If $N=U\g$ is the universal envelope of a Lie algebra then we have
a natural isomorphism $\Alg(U\g,A)\simeq \Lie(\g,Lie(A))$ [CE],
where $\Lie(\g,\h)$ denotes the group of Lie algebra maps $\g\to
\h$ (with pointwise addition) and $Lie(A)$ denotes the  underlying
Lie algebra of the algebra A, that is the Lie algebra where the
Lie bracket is given by $[x,y]=xy-yx$. In our case $Lie (A)$ is an
abelian Lie algebra (since the algebra $A$ is commutative) and
hence $\Lie(\g,Lie(A))\simeq \Lie(\g/[\g,\g],Lie(A))\simeq
\Vec((\g/[\g,\g])^+,A^+)$. Here $\g/[\g,\g]$ is the abelianization
of the Lie algebra $\g$, $\_^+$ is the underlying vector space
functor and $\Vec(V,W)$ is the abelian group of linear maps from
$V$ to $W$ (with pointwise addition). Hence we have:
\begin{Theorem}
Let $G$ be a finite group, and $\bf g$ a Lie algebra. Then
$$
\H^i_{\rm meas}(kG,\Hom(U\g,A))\simeq
\H^i(G,\Vec((\g/[\g,\g])^+,A^+)).
$$
If $|G|^{-1}\in k$ then $\H^i_{\rm meas}(kG,\Hom(U\g,A))=0$.
\end{Theorem}

{\bf\noindent Proof. }The first equality was already explained in
the paragraph preceding the theorem. Note that if $|G|^{-1}\in k$,
then $\Vec((\g/[\g,\g])^+,A^+)$ is uniquely $|G|$ divisible and
hence $\H^i_{\rm meas}(kG,\Hom(U\g,A))=0$. \Qed

\section{Measuring Hopf algebra}

Our aim is to interpret measuring cohomology groups as certain subgroups of
Singer cohomology. We start by showing that there is a way of making the
universal measuring coalgebra $\M(N,A)$ into a Hopf algebra, if $N$ is a Hopf
algebra and $A$ is commutative.

\begin{Proposition}
Let $N$ be a Hopf algebra, $A$ a commutative algebra and $\theta\colon \M(N,A)\ot N\to A$
the universal measuring. There is a unique algebra structure $(\M(N,A),\m,\eta)$
so that $(N,\theta)$ measures $\M(N,A)$ to $A$, i.e. $\theta(\m(f\ot g)\ot n)=
\theta(f\ot n_1)\theta(g\ot n_2)$ and $\theta(\eta\ot\id)=\eta_A\ot\ep_N.$
Moreover $\M(N,A)$ becomes a Hopf algebra with this additional structure.
\end{Proposition}

{\bf\noindent Proof. } Note that the map $\om\colon \M(N,A)\ot
\M(N,A)\ot N\to A$, defined by $\om(f\ot g\ot n)=\theta(f\ot
n_1)\theta(g\ot n_2)$ is a measuring and hence there exists the
unique coalgebra map $\m\colon \M(N,A)\ot \M(N,A)\to \M(N,A)$, so
that $\theta(\m\ot\id)=\om$.

Similarly we define the unit $\eta\colon k\to \M(N,A)$ to
be the unique coalgebra map s.t. $\theta(\eta\ot\id)=\eta_A\ep_N\tau$, where $\tau\colon
k\ot N\to N$ is the natural isomorphism ($x\ot n\mapsto xn$).

A routine computation shows that $\eta$ defined above is a unit for multiplication
$\m$. The following paragraph proves that $\m$ is associative.

Define a measuring $\om_3\colon \M(N,A)\ot\M(N,A)\ot\M(N,A)\ot
N\to A$ by the rule $\om_3(f\ot g\ot h\ot n)=\theta(f\ot n_1)\theta(g\ot
n_2)\theta(h\ot n_3)$ and note that $\theta(\m(\id\ot\m)\ot\id)=\om_3=
\theta(\m(\m\ot\id)\ot\id)$ and hence by the uniqueness $\m(\id\ot
m)=\m(\m\ot\id)$.

Since $\eta$ and $\m$ are coalgebra maps, $\M(N,A)$ is a
bialgebra. Now we conclude the proof by showing the existence of
the antipode. We define $S\colon \M(N,A)^{co-op}\to\M(N,A)$, to be
the unique coalgebra map s.t. $\theta(S\ot\id)=\theta(\id\ot S_N)$
(here $co-op$ refers to the opposite coalgebra structure, i.e.
$\De^{co-op}(f)=f_2\ot f_1$). We claim that $S$ is the antipode.
It is sufficient to show $\theta(S*\id\ot\id)=\theta(\eta\ep\ot
\id)=\theta(\id*S\ot \id)$. This is observed by the following
computation: $\theta(S*\id\ot \id)(f\ot n)=\theta(S(f_1)f_2\ot n)=
\om(S(f_1)\ot f_2\ot n)=\theta(S(f_1)\ot n_1)\theta(f_2\ot n_2)=
\theta(f_1\ot S(n_1))\theta(f_2\ot n_2)=\theta(f\ot S(n_1)n_2)=
\theta(f\ot \ep(n)1_N)=\ep(f)\ep(n)1_A=\theta(\eta\ep(f)\ot n)$;
symmetrically for $\theta(\eta\ep\ot\id)=\theta(\id*S\ot\id)$.
\Qed

\begin{Proposition}
If $N$ is cocommutative, then $\M(N,A)$ is commutative.
\end{Proposition}

{\noindent\bf Proof.} $\th(\m\ot\id)(f\ot g\ot n)=\th(f\ot
n_1)\th(g\ot n_2)=\th(f\ot n_2)\th(g\ot n_1)=\th(gf\ot
n)=\th((\sigma\m)\ot\id)(f\ot g\ot n)$. \Qed

\section{Measuring cohomology vs. Singer cohomology}

In this section we interpret $\H^n_{\rm meas}(T,\Hom(N,A))$ as
$\H^n_n(T,\M(N,A))$.

\begin{Proposition}
If $N$ is left $T$-module bialgebra via $\mu\colon T\ot N\to N$
($t\ot n\mapsto t(n)$), then $\M(N,A)$ is a right $T$-module bialgebra
via $\ol{\mu}\colon \M(N,A)\ot T\to \M(N,A)$ ($f\ot t\mapsto f^t$), which is the
unique map, such that $\th(\ol{\mu}\ot \id)=\th(\id\ot
\mu).$
\end{Proposition}

{\noindent\bf Proof.} By the universal property $\ol{\mu}$ is a
coalgebra map. The following computation proves that $\M(N,A)$ is
also a $T$-module algebra, i.e.
$\ol{\mu}(\m\ot\id)=\m(\ol{\mu}\ot\ol{\mu})\sigma_{2,3}(\id\ot\id\ot\De)$:
$\th((fg)^t\ot n)=\th(fg\ot t(n))=\th(f\ot t_1(n_1))\th(g\ot
t_2(n_2))=\th(f^{t_1}\ot n_1)\th(g^{t_2}\ot
n_2)=\th(f^{t_1}g^{t_2}\ot n)$. \Qed

{\bf\noindent Remark. }{\sl Note that $(T,\M(N,A),\ol\mu,\rho)$,
where $\ol\mu$ is as above and $\rho$ is the trivial coaction
$\rho=\eta\ot\id\colon T\to \M(N,A)\ot T$ is a Singer pair.}

From now on assume $T$ and $N$ are both cocommutative. In this
case we can talk about the differentials $\Reg(T^{\ot
p},\M(N,A))\stackrel{\ol\de}{\to} \Reg(T^{\ot p+1},\M(N,A))$ for computing Sweedler
cohomology $\H^p(T,U(\M(N,A)))$ (here $U\colon {\bf Hopf}
\to {\bf Alg}$ denotes the underlying algebra functor). The following lemma
compares them to Sweedler differentials
$\de\colon
\Reg(T^{n}\ot N,A)\to \Reg(T^{n+1}\ot N,A)$ (recall that we are
identifying $\Hom(T^{i}\ot N,A)\simeq \Hom(T^{i},\Hom(N,A))$.

\begin{Lemma}
If $\ol{\alpha}\colon T^{\ot n}\to \M(N,A)$ is the coalgebra map
corresponding to the measuring $\alpha\colon T^{\ot n}\ot N\to A$
then, $\ol\de\ol{\al}=\ol{\de\al}$ i.e. the unique coalgebra map
corresponding to the measuring $\de\alpha\colon T^{\ot n+1}\ot
N\to A$ is $\ol{\de}\ol{\al}\colon T^{\ot n+1}\to \M(N,A)$.
\end{Lemma}

{\noindent\bf Proof. } By the universal property it is sufficient to prove that
$\th(\ol{\de}\ol{\alpha}\ot\id)=\de\alpha$. Since
$\th(\ol\al*\ol\beta\ot\id)=\ol{\al*\beta}$ and also
$\ol{\al^{-1}}=\ol{\al}(\id_{T^{\ot p}}\ot S)=S\ol\al={\ol{\al}}^{-1}$
(it is a coalgebra map, since $T$ is cocommutative), it is
enough to see that $\theta(\ol{d_i}\ol{\alpha}\ot\id_N)=d_i \al$, where
$\ol{d_i}\ol{\alpha}=\ol\alpha(\id_T^{i}\ot\m_T\ot\id_T^{n-i-1})\colon T^{\ot n+1}\to \M(N,A)$,
$d_i\al=\al(\id_T^{i}\ot\m_T\ot\id_T^{n-i-1}\ot\id_N)\colon T^{\ot n+1}\ot N\to A$,
for $i=0,\ldots, n$ and $\ol{d_{n+1}}\ol{\alpha}=\ol{\mu}(\ol\al\ot\id_T)$,
$d_{n+1}\al=\al(\id_T^{n}\ot\mu)$. Let $D_i=\id^{i}_T\ot\m_T\ot\id_T^{n-i-1}\ot\id_N\colon
T^{\ot n+1}\ot N\to T^{\ot n}\ot N$ for $i=0,\ldots, n$ and $D_{n+1}=\id_T^n\ot\mu$ and note that
$d_i\al=\al D_i$ and hence $\theta(\ol\al\ot \id_N)D_i=d_i\al$. \Qed

Observe that $\ol{\eta_A\ep_{T^{\ot p}\ot
N}}=\eta_{\M(N,A)}\ep_{T^{\ot p}}$ and hence

$$\ol{(\_)}\colon
\{\alpha\colon T^{\ot p}\ot N\to A|\al\; {\rm measures}\}\to
\Coalg(T^{\ot p},\M(N,A))$$ gives an isomorphism of complexes
$$(\Reg_{meas}(T^{\ot p},\Hom(N,A)),\de)\to (\Coalg(T^{\ot
p},\M(N,A)),\ol\de).$$

\begin{Theorem}
$\H^n_n(T,\M(N,A))\simeq \H^n_{\rm meas}(T,\Hom(N,A))$ for $n\ge 2$.
\end{Theorem}
{\noindent\bf Proof. }Apply Proposition \ref{desc} and the lemma above. \Qed

{\bf\noindent Remark. }{\sl Hence the degree two measuring
cohomology characterizes those Hopf algebra extensions $\M(N,A)\to
H\to T$, for which there exists an $\M(N,A)$-module coalgebra map
$\xi\colon H\to \M(N,A)$. It is also easy to see that
$\Z^1_1(T,\M(N,A))=\Z^1(T,\M(N,A))=\Z^1_{\rm meas}(T,\Hom(N,A))$ and
that $B^1_1(T,\M(N,A))=\B^1(T,\M(N,A))=0$. Hence $\H^1_{\rm meas}(T,\Hom(N,A))$ is a quotient of
$\H^1_1(T,\M(N,A))=
\H^1(T,\M(N,A))$.}

\section{Five term exact sequence for a smash product}
The purpose of this section is to prove the following theorem by
explicitly describing the maps involved.

\begin{Theorem}
Let $H=N\rt T$ be a smash product of cocommutative Hopf algebras
(more precisely, we are given an action $\mu\colon T\ot N\to N$,
that makes $N$ into a $T$-module bialgebra) and let the commutative
algebra $A$ be a trivial $H$-module. Then we have the following
exact sequence:
\begin{eqnarray*}
0&\to& \H^1_{\rm meas}(T,\Hom(N,A))\stackrel{\io}{\to} \tilde{\rm
H}^2(H,A)\stackrel{\res}{\to}
\H^2(N,A)^T\\
&\stackrel{d}{\to}& \H^2_{\rm meas}(T,\Hom(N,A))\stackrel{j}{\to}
\tilde{\rm H}^3(H,A).
\end{eqnarray*}
\end{Theorem}

We prove the above theorem  by transporting some arguments from
[Ta] into our more general setting.

First we have to define the Hopf algebra analog of the stable part of
cohomology.

\begin{Definition}
Let $N$, $T$, $\mu$, $A$ be as above. We say that a cohomology
class $[f]\in \H^i(N,A)$, where $f\in \Z^i(N,A)$, is $T$-stable if
there exists a convolution invertible linear map $g\colon T\ot
N^{\ot i-1}\to A$, such that $f*(f^{-1})^t=\de' g(t\ot\_)$. Here
$\de'$ denotes the Sweedler differential (see Section \ref{s2})
from $\Reg(T\ot N^{\ot i-1},A)\simeq \Reg(N^{\ot i-1},\Hom(T,A))$ to
$\Reg(T\ot N^{\ot i},A)\simeq \Reg(N^{\ot i},\Hom(T,A))$. The
subgroup of $\H^i(N,A)$ consisting of all $T$-stable elements is
called the $T$-stable part of cohomology and is denoted by
$\H^i(N,A)^T$.
\end{Definition}
Note that if $T=kG$ then $\H^i(N,A)^T=\H^i(N,A)^G$.

The following lemma is the main tool in establishing this result.
It is a generalization of the essential part of Proposition 1 from
[Ta]:
\begin{Lemma}\label{L41}
Let the Hopf algebra $H$ be a smash product of cocommutative Hopf
algebras $N$ and $T$. Furthermore assume that $H$ acts trivially
on a commutative algebra $A$. Then every cocycle $f\colon H\ot
H\to A$ is cohomologous to a cocycle $f'\colon H\ot H\to A$, which
is trivial on $N\ot T$, i.e. $f'(n\ot t)=\ep(n)\ep(t)1_A$.
\end{Lemma}

{\bf\noindent Proof }(of Lemma \ref{L41}). Let the extension
$$A\stackrel{\io}{\to} K
\displaystyle{\mathop{\rightleftarrows}^{\pi}_\chi} H$$
\noindent be an $H$-comodule algebra extension with associated
$2$-cocycle $f$ (see Section \ref{s1}). We shall denote the
$H$-comodule structure on $K$ by $\rho\colon K\to K\ot H$. We will
show that $\chi$ can be \lq\lq repaired" into a $\chi'\colon H\to K$
that satisfies the equality $\chi'(nt)=\chi'(n)\chi'(t)$, for
$n\in N$, $t\in T$. Then it is easy to see that the cocycle
$f'=f'_{\chi'}\colon H\ot H\to A$ associated to $\chi'$ satisfies the desired
condition.

Let $\{u_i\}_{i\in I}$ be a basis for $N$ and let $\{v_j\}_{j\in
J}$ be a basis for $T$. Then $\{u_i\rt v_j\}_{(i,j)\in I\times J}$
is a basis for $H$. Define a linear map $\chi'\colon H\to K$ by
the rule $\chi'(u_i\rt v_j)=\m_K(\chi(u_i)\ot\chi(v_j))$. The
following calculation shows that $\chi'$ is an $H$-comodule
map, i.e. $\rho\chi'=(\chi'\ot\id)\De_H$:
\begin{eqnarray*}
(\chi'\ot\id)\De(u_i\rt v_j)&=& \sum \chi'((u_i)_1\rt (v_j)_1)\ot
(u_i)_2\rt (v_j)_2 \\
&=& \sum \chi(((u_i)_1)\chi((v_j)_1)\ot (u_i)_2 (v_j)_2 \\
&=& \sum \m_{K\ot H}((\chi((u_i)_1)\ot
(u_i)_2)\ot(\chi((v_j)_1)\ot (v_j)_2))\\
&=& \m_{K\ot H}(
\rho\chi(u_i)\ot\rho\chi(v_j)) =
\rho\m_K (\chi(u_i)\ot \chi(v_j))\\
&=& \rho\chi'(u_i\rt v_j).
\end{eqnarray*}
Now observe $\chi'(u_i)=\chi(u_i)$ and $\chi'(v_j)=\chi(v_j)$ and
hence $\chi'(nt)=\chi'(n\rt t)=\chi'(n)\chi'(t)$ for $n\in N$ and
$t\in T$. \Qed

A cocycle $f'$ that satisfies the condition of Lemma \ref{L41}
will be called a normalized cocycle.

\begin{Corollary} \label{C1}
Let $f\colon H\ot H\to A$ be a normalized cocycle, where $H$ is a
smash product of $N$ and $T$ acting trivially on the commutative
algebra $A$. Then $f$ satisfies the following equations:
\begin{eqnarray}
f(nt\ot h') &=& \sum f(t_1\ot h'_1)f(n\ot t_2 h'_2)\label{ee1}\\
f(nt\ot t') &=& \ep(n)f(t\ot t')\label{ee2}\\
f(h\ot n't') &=& \sum f(h_1\ot n'_1)f(h_2 n'_2\ot t')\label{ee3}\\
f(n\ot n't') &=& f(n\ot n')\ep(t')\label{ee4}\\
f(nt\ot n't') &=& \sum f(t_1\ot t')f(t_2\ot n'_1)f(n\ot t_3(n'_2))\label{ee5}
\end{eqnarray}
for $n,n'\in N$, $t,t'\in T$ and $h,h'\in H$.
\end{Corollary}

{\bf\noindent Proof. }The equations (\ref{ee1}) and (\ref{ee3}) are just special
cases of the cocycle condition. Equations (\ref{ee2}) are (\ref{ee4}) are special
cases of (\ref{ee1}) and (\ref{ee3}) respectively and (\ref{ee5}) follows from (\ref{ee1})-(\ref{ee4}).
\Qed

\begin{Corollary} \label{C2}
A map $f\colon H\ot H\to A$ is a normalized $2$-cocycle if and
only if the following are satisfied:
\begin{enumerate}
\item
$f|_{N\ot T}=\ep$
\item
$f|_{N\ot N}$ is a $2$-cocycle on $N$
\item
$f|_{T\ot T}$ is a $2$-cocycle on $T$
\item
$f(tt'\ot n')=\sum f(t'_1\ot n'_1)f(t\ot t'_2(n'_2))$, where
$n'\in N$ and $t,t'\in T$
\item
$\sum f(n_1\ot n'_1)f^{-1}(t_1(n_2)\ot t_2(n'_2))=\sum f(t_1\ot
n_1)f^{-1}(t_2\ot n_2 n'_1)f(t_3\ot n'_2)$, where $n,n'\in N$ and
$t\in T$.
\end{enumerate}
Moreover, the data $f|_{N\ot N}$, $f|_{T\ot T}$, $f|_{T\ot N} $
satisfying the conditions above determine a unique normalized
cocycle.
\end{Corollary}

{\bf\noindent Proof.} First assume that $f$ is a normalized
$2$-cocycle on $H$. Then clearly $f$ is also a $2$-cocycle on both
$N$ and $T$. The equations in conditions \ref{ee3} and \ref{ee5} are obtained from
the cocycle condition together with the equations (\ref{ee2}) and (\ref{ee3}) from
the previous corollary.

Now suppose that we have the data from conditions \ref{ee1}-\ref{ee5}. Then
the Equation (\ref{ee5}) of the Corollary \ref{C1}
gives a formula for a map $f\colon H\ot H\to A$ (it is well defined, since it is linear
in each of the variables). An elementary (but lengthy) computation
shows that the cocycle condition is satisfied. \Qed \vskip 1em

Let $\tilde{\H}^i(H,A)$ be the kernel of the restriction
homomorphism  $\H^i(H,A)\stackrel{\res}{\to} \H^i(T,A)$. Since the
inclusion $T\to H$ splits we can conclude that
\begin{Proposition}
$$\H^i(H,A)\simeq
\H^i(T,A)\oplus \tilde{\H}^i(H,A).$$
\end{Proposition}
We shall denote the group of
normalized cocycles $H\ot H\to A$ that are trivial when restricted
to $T$ by ${\rm Z'}^2(H,A)$, i.e. ${\rm Z'}^2(H,A) = \{f\in
\Z^2(H,A)| f(n\ot t)=\ep(n)\ep(t)\;{\rm and}\; f(t\ot
t')=\ep(t)\ep(t'),\;n\in N,\; t,t'\in T\}$. Furthermore, let ${\rm
B'}^2(H,A)=\B^2(H,A)\cap \Z'^2(H,A)$. Using the canonical map
$H\to T$ together with Corollaries \ref{C1} and \ref{C2} we can
show that there is an injective map $\B^2(T,A)\to \B^2(H,A)$ and
hence
\begin{Proposition}
$\H'^2(H,A)=\Z'^2(H,A)/\B'^2(H,A)\simeq \tilde{\rm H}^2(H,A)$.
\end{Proposition}

We proceed by defining the homomorphisms involved in the
generalized Tahara sequence, and also prove exactness at the same
time. \vskip 1em

{\noindent\bf The injective homomorphism $\H^1_{\rm
meas}(T,\Hom(N,A))\stackrel{\io}{\to}\tilde{\rm H}^2(H,A)$:}

Define a homomorphism $\io\colon \Z^1_{\rm
meas}(H,\Hom(N,A))\to \Z'^2(H,A)$ by the rule $\io(f)(nt\ot
n't')=f(t)(n')\ep(n)\ep(t')$. This homomorphism induces an injective
homomorphism $\H^1_{\rm meas}(T,\Hom(N,A))\to \tilde{\rm
H}^2(H,A)$. \Qed \vskip 1em

{\noindent\bf The exactness at $\tilde{\rm
H}^2(H,A)$:}\nopagebreak

We claim, that the image of the homomorphism just described equals
the kernel of the restriction homomorphism $\tilde{\rm
H}^2(H,A)\stackrel{\res}{\to}\H^2(N,A)^T$.

Clearly $\res\io=0$. Suppose the cocycle $f\in \Z'^2(H,A)$ is such
that $f|_{N\ot N}\in \B^2(N,A)$, that is there exists $g\colon
N\to A$ such that $f(n\ot n')=\delta g(n\ot n')$. Extend $g$ to a
linear map $g\colon H\to A$ by the rule $g(n\rt t)=g(n)\ep(t)$. Now
define $f'\in \Z^2_{\rm meas}(T,\Hom(N,A))$ by the rule
$f'(t)(n')=\sum f(t_1\ot n'_1)g^{-1}(n'_2)g(t_2(n'_3))$ A routine
calculation shows that $f*\delta g^{-1}=\io(f')$ (we use
Equation (5) of Corollary \ref{C2} to expand $f(nt\ot n't')$ and
take into account that $f|_{N\ot N}=\delta g$ and that $f|_{T\ot
T}=\ep$) and hence $[f]=\io([f'])$. \Qed \vskip 1em

{\noindent\bf The homomorphism $\H^2(N, A)^T\stackrel{d}{\to}
\H^2_{\rm meas}(T,\Hom(N,A))$ and the exactness at $\H^2_{\rm
meas}(T,\Hom(N,A))$:}

Take $[f]\in H^2(N, A)^T$, for $f\in \Z^2(N,A)$. Then there is a
$g\colon T\ot N\to A$ s.t. $f*(f^{-1})^{t}=\de g(t\ot\_)$, i.e.
$\sum f(n_1\ot n'_1)f^{-1}(t_1(n_2)\ot t_2(n'_2))=\sum g(t_1\ot
n_1)g^{-1}(t_2\ot n_2 n'_1)g(t_3\ot n'_2)$. Now define $d$ by
$(df)(t\ot t'\ot n)=\sum g(t'_1\ot n_1)g^{-1}(t_1t'_2\ot
n_2)g(t_2\ot t'_3(n_3))$. Clearly $df\in \Z^2_{\rm
meas}(T,\Hom(N,A))$. We claim that the class $[df]$ is independent
of the choice of $g$, which in turn also implies that
$d(\B^2(N,A))\subseteq \B^2(T, \Hom(N,A))$ and hence $d$ gives
rise to a homomorphism $\H^2(N, A)^T\to \H^2_{\rm
meas}(T,\Hom(N,A))$. So suppose there is a $g'\colon T\ot N\to A$
such that $\sum f(n_1\ot n'_1)f^{-1}(t_1(n_2)\ot t_2(n'_2))=\sum
g'(t_1\ot n_1)g'^{-1}(t_2\ot n_2n'_1)g'(t_3\ot n'_2)$. We need to
show that there exists $w\in \Reg_{\rm meas}(T,\Hom(N,A))$ such
that $\sum g(t'_1\ot n_1)g^{-1}(t_1t'_2\ot n_2)g(t_2\ot
t'_3(n_3))g'^{-1}(t'_4\ot n_4)g'(t_4t'_5\ot n_5)g'^{-1}(t_5\ot
t'_6(n_6)) = \sum w(t'_1\ot n_1) w^{-1}(t_1t'_2\ot n_2)
w(t_2\ot t'_3(n_3))$. Define $w$ by $w(t\ot n)=\sum g(t_1\ot
n_1)g'^{-1}(t_2\ot n_2)$ and observe that it does the trick.

It is clear that $d\io=0$. Suppose $df\in \B^2(T,\Hom(N,A))$. Then
there exists a $w\in \Reg_{\rm meas}(T,\Hom(N,A))$ s.t. $(df)(t\ot
t'\ot n)=\sum w(t'_1\ot n_1)w^{-1}(t_1t'_2\ot n_2)w(t_2\ot
t'_3(n_3))$. Define $z\colon T\ot N\to A$ by
$z(t\ot n)=\sum h(t_1\ot n_1)u^{-1}(t_2\ot
n_2)$ and note that it gives rise to a normalized cocycle $z\in
\Z'^2(H,A)$ given by $z(nt\ot n't')=z(t\ot n)$ (see Lemma
\ref{C2}) and that $i[z]=[f]$. \Qed\vskip 1em

{\noindent\bf The homomorphism $\H^2_{\rm meas}(T,\Hom(N,A))\stackrel{j}{\to}
\tilde{\rm H}^3(H,A)$ and the exactness at $H^2_{\rm
meas}(T,\Hom(N,A))$}

Let $f\in \Z^2_{\rm  meas}(T,\Hom(N,A))$. Define a map $jf\colon
H\ot H\ot H\to A$ by $jf(nt\ot n't'\ot n''t'')=f(t\ot t'\ot n'')$.
A straightforward calculation shows that $jf$ is a $3$-cocycle on
$H$. Suppose that $f$ is a measuring $2$-coboundary. Then there
exists $v\in \Reg_{\rm meas}(T,\Hom(N,A))$ s.t. $f(t\ot t'\ot
n'')= \sum v(t'_1\ot n''_1)v^{-1}(t_1t'_2\ot n''_2)v(t_2\ot
t'_3(n''_3))$. Now let $v'\colon H\ot H\to A$ be defined by
$v'(nt\ot n't')=v(t\ot n')$ and show $jf=\delta^2 v'$. This proves
that the homomorphism $j$ is well defined.

Suppose $[h]\in \H^2(N,A)^T$ and let $u\colon T\ot N\to A$ be s.t.
$h*(h^{-1})^t=\delta'(u(t\ot\_))$. Define $v\colon H\ot H\to A$ by
$v(tn\ot t'n')=\sum u(t_1\ot n'_1)h(n\ot t_2(n'_2))$ and observe
that $jdh=\delta^2 v\in B^3(H,A)$. This shows that $jd=0$.

Now suppose the measuring $2$-cocycle $f$ is such that $jf$ is a
$3$-coboundary. Then there exists $v\in \Reg(H\ot H,A)$ s.t.
$jf=\delta^2 v$. Define a $u\colon T\ot N\to A$ by $u(t\ot n)=\sum
v(t_1\ot n_1)v^{-1}(t_2(n_2)\ot t_3)$ and $h\colon N\ot N\to A$ by
$h(n\ot n')=v(n\ot n')$ and note that $\delta^1(h)=\eta\ep$ (hence
$h\in \Z^2(N,A)$) and that $h*(h^{-1})^t=\de' u(t\ot\_)$, so that
$[h]\in \H^2(N,A)^T$. Observe also that $f(t\ot t')= \sum
u(t'_1)u^{-1}(t_1t'_2)u^{t'_3}(t_2)$, i.e $[f]=d[h].\;\Qed$

\begin{Corollary}
If $H=N\ot T$, i.e. if the action of $T$ on $N$ is trivial, then
we have a canonical isomorphism
$$
\H^2(H,A)\simeq \H^2(T,A)\oplus \H^1_{\rm meas}(T,\Hom(N,A))\oplus
\H^2(N,A).
$$
\end{Corollary}

\begin{Proposition}
If the action of $T$ on $N$ is trivial, then
$$
\H^1_{\rm meas}(T,\Hom(N,A))\simeq \Pp(T,N,A),
$$
where $\Pp(T,N,A)$ denotes the abelian group of maps $f\colon T\ot
N\to A$ that measure in both variables, i.e correspond to algebra
maps $T\to \Hom(N,A)$ and $N\to \Hom(T,A)$.
\end{Proposition}

{\bf\noindent Remark.}{\sl The isomorphism $\H^2(N\ot T, A)
\stackrel{\sim}{\to} \H^2(T,A)\oplus \H^2(N,A)\oplus \Pp(T,N,A)$
has a description similar to that in the case of group cohomology
(for the group cohomology case see for instance [Kar]):
$[f]\mapsto ([f|_{T\ot T}],[f|_{N\ot N}],\tilde{f})$,  where
$\tilde{f}(t\ot n)=\sum f(t_1\ot n_1)f^{-1}(n_2\ot t_2)$.}

\section{On $\H^2(kG\rtimes U{\bf g},A)$}

Here we illustrate how the generalized Tahara sequence sheds some
light on the Sweedler cohomology, when the cocommutative Hopf
algebra in question is a smash product of a group algebra $T=kG$
and the universal envelope of a Lie algebra $N=U\g$. This is
always the case when $k$ is algebraically closed and of
characteristic $0$ [Gr].

\begin{Theorem}\label{sdLg}
Let $G$ be a finite group acting on a Lie algebra ${\bf g}$,
furthermore let $A$ be a commutative algebra which is also a
trivial $U{\bf g}\rt kG$ module and assume the ground field $k$
contains $|G|^{-1}$. Then
$$
\H^2(U{\bf g}\rt kG)\simeq \H^2({\bf g},A^+)^G\oplus \H^2(G,A).
$$
\end{Theorem}

{\bf\noindent Proof. } If $|G|^{-1}\in k$ then $\H^i_{\rm
meas}(kG,\Hom(U\g,A))$ is trivial and hence the restriction
homomorphism $\res\colon \tilde{\rm H}^2(N\rt T)\to \H^2(N,A)^T$
is an isomorphism. So we get the isomorphism $\H^2(N\rt T,A)\simeq
\H^2(T,A)\oplus \H^2(N,A)^T$. Now $\H^2(T,A)=\H^2(G,\U(A))$,
$\H^2(N,A)\simeq \H^2(\g,A^+)$ (see [Sw1]) and $\H^2(N,A)^T\simeq
\H^2(\g,A^+)^G.\;\Qed$

\begin{Example} Assume the ground field $k$ has characteristic $0$,
let ${\bf g}={sl}_n(k)$ be a Lie algebra consisting of trace zero
$n\times n$ matrices (with Lie bracket given by a commutator) and
let $G\simeq C_n\le Gl_n(k)$ be a group generated by the standard
$n$-cycle permutation matrix acting on ${\bf g}$ by conjugation.
Then $\H^2(U{\bf g}\rt k G, k)\simeq \kb/(\kb)^n$
\end{Example}

{\bf\noindent Proof. }Apply Theorem \ref{sdLg} and note that
$\H^2({\bf g},A^+)$ is trivial by the Whitehead's second lemma
[We] and that $\H^2(C_n,\kb)=\kb/(\kb)^n$ [We]. \Qed

\end{document}